\documentclass[11pt]{article}
\usepackage[cp1251]{inputenc}
\usepackage[T2A]{fontenc}
\usepackage[russian,english]{babel}
\usepackage{amsfonts}
\usepackage{amsmath}
\usepackage{amssymb}
\usepackage{graphicx}
\newcounter{assumption}
\newcounter{definition}
\newcounter{lemma}
\newtheorem{theorem}{Theorem}

\newtheorem{assumption}{Assumption}
\newtheorem{corollary}{Corollary}
\newtheorem{definition}{Definition}

\newtheorem{lemma}{Lemma}

\newenvironment{proof}[1][Proof]{\textbf{#1.} }{\ \rule{0.5em}{0.5em}}

\newcommand{\R}{{\mathbb R}}
\renewcommand{\C}{\mathbb C}
\newcommand{\Z}{{\mathbb Z}}
\newcommand{\N}{{\mathbb N}}

\newcommand{\eps}{{\varepsilon}}

\title{Nonautonomous gradient-like ODEs on the  circle:\\
classification, structural stability and autonomization}
\author{L.M. Lerman, E.V. Gubina\\
\normalsize
Institute of Information Technology, Mathematics and Mechanics,\\
\normalsize
Lobachevsky National Research State University of Nizhny Novgorod, \\
\normalsize e-mail: lermanl@mm.unn.ru \\
}
\date{}
\begin{document}
\maketitle
\vspace{3.mm}
\begin{flushright}
To J\"urgen, with thanks and best wishes
\vspace{3.mm}
\end{flushright}

\begin{abstract} We study a class of scalar differential equations on the
circle $S^1$. This class is characterized mainly by the property that any
solution of such an equation possesses exponential dichotomy both on the
semi-axes $\R_+$ and $\R_+$. Also we impose some other assumptions on the
structure of the foliation into integral curves for such the equation.
Differential equations of this class are called gradient-like ones.
As a result, we describe the global behavior of the foliation, introduce
a complete invariant of uniform equivalency, give standard models for
the equations of the  distinguished class. The case of almost periodic
gradient-like equations is also studied, their classification is presented.
\end{abstract}

\begin{flushleft}

\vspace{3.mm}
\end{flushleft}

\section{Introduction}

In 1973 a short note \cite{LeSh} was published where for nonautonomous vector fields (NVFs)
given on a smooth closed manifold $M$
a definition of uniform equivalency of two such NVFs was given and on this basis the structural
stability of nonautonomous vector
fields was defined. When $\dim M = 2$ a class of structurally stable nonautonomous vector fields was selected,
the invariant determining
the uniform equivalence was found, Morse type inequalities connecting the topology of $M$ and
the set of its integral curves were derived.
When these results were announced, articles on nonautonomous nonlinear dynamics were rather rare
(except, of course, periodic in time
systems), though they existed, see, for instance, \cite{Markus,Sell}.
Now the interest to the nonautonomous
dynamics became much more intensive, several recent books demonstrate this
\cite{KR,Cheban,Potzsche}. This gave us a kick to return to this theme.

In this paper we present proofs of some statements of \cite{LeSh} for the one dimensional situation,
that is when $M$ is the circle $S^1$. It is worth noting that one dimensional case itself was not
discussed in \cite{LeSh} at all.

We study a class of nonautonomous scalar differential equations of the type
$$\dot x = f(t,x),$$ where scalar function $f$ is 1-periodic in $x$ and defines a uniformly continuous map
$\R \to V^r(S^1)$
to the space of $C^r$–smooth 1-periodic functions in $x$, $r\ge 1$. We consider this space with
$C^r$-norm, $r\ge 1,$ thus $V^r(S^1)$ becomes the Banach space. The space of such maps $\R\to V^r(S^1)$ becomes
a Banach space if we endow it with the norm ${\bf |}v{\bf |}=\sup\limits_{t\in\R}||v(t)||$
where $||\cdot ||$ means the norm in $V^r(S^1).$  An ordinary differential equation (ODE)
determines a foliation of the extended phase space $\R \times S^1$ into integral curves (graphs of the solution).
Henceforth, we consider the manifold $\R \times S^1$ with its uniform
structure of the direct product of the standard uniform structure on $\R$ given by the metrics $|x-y|$
and metrics on $S^1 = \R/\Z$ induced from $\R.$

Two such ODEs are called {\em uniformly equivalent} if their related foliations are equimorphic (see below
in more details). Below we
shall present four assumptions on the class of ODEs under consideration. When they hold
we call the equation to be gradient-like one. For equations from
this class we prove the following properties:
\begin{itemize}
\item A gradient-like ODE is structurally stable w.r.t. perturbations within
the class of uniformly continuous maps $\R \to V^r(S^1)$;
\item ODEs of this class are rough that is a conjugating equimorphism can be chosen close to the identity
map $id_{\R\times S^1}$ if a perturbation is small enough;
\item A combinatorial type invariant is introduced which is the complete
invariant of the uniform conjugacy for the gradient-like ODEs;
\item for each gradient-like ODE there is an uniformly equivalent
asymptotically autonomous ODE on $S^1$;
\item if a gradient-like ODE is in addition almost periodic in $t$
uniformly w.r.t. $x$, then this ODE is uniformly equivalent to an
autonomous scalar ODE with simple zeroes.
\end{itemize}
The results can be also formulated and proved similarly for the case $M=I$, a segment, as well,
in this case it is assumed the boundary curves $\mathbb R\times \partial I$ be integral curves
of the differential equation to avoid some complications.
The exposition in the paper is carried out for the case $S^1$.

The beginning of the study of nonautonomous vector fields from the viewpoint of their roughness and their
structure was initiated by L. Shilnikov, our joint note \cite{LeSh} was the first result in this direction.
The idea itself on the necessity of extending the notion of roughness (structural
stability) \cite{AndP} onto nonautonomous vector fields goes back to A.A. Andronov and was publicized by his wife and
collaborator E.A. Leontovich-Andronova in her talk at the III All-Union USSR Mathematical Congress \cite{EA}.

The first problem here was to find a proper equivalency relation for two nonautonomous vector fields
which could be served for a basis
of the classification. By that time there existed mainly two approaches to the study of nonautonomous
systems. The first one goes back to
Bebutov and his translation dynamical system in the space of bounded continuous functions
(see \cite{NemS}).
That approach was transformed later into the theory of skew product systems \cite{Sell}.
The second approach was the study of a nonautonomous system itself but the problem
of the equivalency of two such systems was not formulated that time.

\section{Nonautonomous vector fields}

Let $M$ be a $C^{\infty }$-smooth closed manifold and ${\cal
V}^{r}(M)$ be the Banach space of $C^r$-smooth vector fields on
$M$ endowed with $C^r$-norm. A \textit{$C^r$-smooth nonautonomous
vector field} on $M$ (below NVF, for brevity) is an uniformly continuous bounded map $v:
{\mathbb R} \to {\cal V}^{r}(M)$. If this map $v$ is also $C^s$-differentiable map, whose derivatives
up to order $s$ are uniformly
continuous, we call $v$ to be a $C^{r,s}$-smooth nonautonomous vector field. Every nonautonomous vector
field $v$ defines its solutions
$x_{t,x_0}: \mathbb R \to M$ being $C^1$ maps, if $s=0$ and $C^{s+1}$ maps, if $s>0$. As a mapping
from $\mathbb R \times \mathbb R\times M \to
M$, $(t;\tau,x_0)\to x(t;\tau,x_0)\in M,$ this mapping is of $C^{min\{r,s\}}$. Solutions of the
vector field generate a foliation
$F$ of the manifold $M\times \R$ (extended phase space) into its integral curves $\cup_t(x(t),t).$
Henceforth we consider the manifold $M\times \mathbb R$
with its standard uniform structure (see, for instance \cite{Kelley}). All uniformly continuous maps
 of $M\times \mathbb R$ to itself
are considered with respect to this uniform structure. Recall that a
homeomorphism $h: M\to M$ of a uniform space is called an equimorphism, if
both mappings $h, h^{-1}$ are uniformly continuous.

An equivalence relation for nonautonomous vector fields proposed in \cite{LeSh}
is as follows.
\begin{definition}
Two NVFs $v_1, v_1$ are \textbf{uniformly equivalent} if foliations $F_1, F_2$ are uniformly
equivalent, that is there is an equimorphism
$h: M\times\R \to M\times \R$ respecting foliations
(i.e. sending every integral curve $\gamma$ of $F_1$ to an integral curve of $F_2$ preserving
its orientation in $\R$).
\end{definition}
It is clear that this equivalency relation distinguishes NVFs in which the asymptotic behavior
of integral curves are different.
This relation allows us to introduce the notion of structurally stable NVFs.
\begin{definition}
An NVF $v$ is called \textbf{structurally stable}, if there is a neighborhood $\mathcal U$ of
$v$ in the space $NV$ such that
all NVFs in this neighborhood are uniformly equivalent.
\end{definition}
The development of the modern theory of (autonomous) dynamical systems showed that the topological
equivalency relation is in fact too rigid to get a classification of multidimensional dynamical
systems: the structure of such a system can be extremely complicated \cite{GST93a,GST99} due,
in particular, to the phenomena like Newhouse ones \cite{N74,N79,GST93b}.

Nonetheless, this equivalency relation is good enough to classify relatively simple (Morse-Smale)
systems, though even in this case the classifying
invariant become comparable in its complexity with the structure of a systems itself \cite{GMP}.
This is the achievements of
the recent time, at times when this theory started, it was not known.

\section{One dimensional gradient-like NVFs}

From now on in this paper we consider NVFs on the unique closed smooth 1-dimensional manifold,
that is on a circle $S^1 = \mathbb R/\mathbb Z$. We formulate first the restrictions on NVFs which will allow one
to classify these vector fields. Let $x$ be a (1-periodic) coordinate on $S^1$,
then a NVF defines a scalar ordinary differential equation (ODE, for brevity)
\begin{equation}\label{1d}
\dot x = f(t,x),\; x\in S^1
\end{equation}
with 1-periodic in $x$ function $f.$ Here $f$ is continuous
in both variables and uniformly continuous w.r.t. $t$, differentiable in $x$
uniformly w.r.t. $t$. Hence, the usual existence and uniqueness
theorem is valid for this ODE and any its solution is extended in $t$ on the whole $\mathbb R$.
Recall that an integral curve of the equation (below IC, for brevity) is
the graph of the related solution in the extended phase space
$\mathbb R\times S^1$. Thus the extended phase space is foliated into integral
curves. We intend to classify these equations w.r.t. the uniform equivalency relation using some combinatorial
invariant to be defined later on. To do this, we impose some restrictions on the class of NVFs under consideration.

The first of these restrictions is the following
\begin{assumption}
Any solution of $(\ref{1d})$ possesses exponential dichotomy both on $\R_+$ and $\R_-$ $\cite{MS}$.
\end{assumption}
Recall the related definitions \cite{MS,Copp,Demid}.
\begin{definition}
Let $x(t)$ be a solution of the equation $(\ref{1d})$. One says that this solution satisfies
the exponential dichotomy of the stable type on the semi-axis $\R_+ = \{t\ge 0\}$, if there are positive
constants $C,\lambda$ such that the linearized at this solution linear ODE
$\dot\xi = a(t)\xi,$ $a(t)=f_x(t,x(t)),$ satisfies the inequality
$$\displaystyle{\exp[\int\limits_{\tau}^{t} a(s)ds] \leq C\exp[-\lambda(t-\tau)]}
$$
for all $t,\tau,\;t\ge \tau\ge 0,$
\end{definition}

\begin{definition}
Let $x(t)$ be a solution of the equation $(\ref{1d})$. One says that this solution satisfies
the exponential dichotomy of the unstable type on the semi-axis $\R_+$, if there are
positive constants $C,\lambda$ such that the linearized at this solution linear ODE
$\dot\xi = a(t)\xi$ satisfies the inequality
$$\displaystyle{\exp[\int\limits^{\tau}_{t} a(s)ds] \leq C\exp[\lambda(t-\tau)]}
$$
for all $t,\tau,\;0\le t\le \tau.$
\end{definition}

Similar notions of exponential dichotomy of the stable and unstable types are defined for the semi-axis $t\le
0$. To be more precise, let us present these definitions as well.
\begin{definition}
Let $x(t)$ be a solution of the equation $(\ref{1d})$. One says that this solution satisfies
the exponential dichotomy of stable type on the semi-axis $\R_-$, if there are positive
constants $C,\lambda$ such that the linearized at this solution equation $\dot\xi = a(t)\xi$
 satisfies the inequality
$$
\displaystyle{\exp[\int\limits_{\tau}^{t} a(s)ds] \leq C\exp[-\lambda(t-\tau)]}
$$
for all $\tau \le t \le 0$.
\end{definition}

\begin{definition}
Let $x(t)$ be a solution of the equation $(\ref{1d})$. One says that this solution satisfies
the exponential dichotomy of the unstable type on the semi-axis $\R_-$, if there are
positive constants $C,\lambda$ such that the linearized equation $\dot\xi = a(t)\xi$
 satisfies the inequality
$$
\displaystyle{\exp[\int\limits^{\tau}_{t} a(s)ds] \leq C\exp[\lambda(t-\tau)]}
$$
for all $t\le \tau \le 0.$
\end{definition}

It follows from the Hadamard-Perron theorem \cite{Anosov} that for a
solution $\gamma(t)$ to $(\ref{1d})$ possessing the exponential dichotomy of stable type on $\R_+$ there is its
uniform neighborhood $U\subset \R_+\times S^1$ such that all solutions of the equation
(\ref{1d}) starting at $t=t_0\ge 0$ within $U$ exponentially fast
tend to this IC as $t\to \infty.$
Moreover, such $U$ can be chosen in such a way that a quadratic Lyapunov function would exist in $U$
and lateral boundary curves of this $U$ are the level lines of the Lyapunov function.
Recall how to construct such a function (see, for instance, \cite{Demid}). Let us write down the equation near
the integral curve using local variables $x = \gamma(t)+u$. Then the equation transforms into the form
\begin{equation}\label{linear}
\dot u = a(t)u + h(t,u),\;a(t)= f_x(t,\gamma(t)),\;h(t,0)=h_u(t,0)\equiv 0.
\end{equation}
Since the map $\hat f:\R \to V^r(S^1),$ $r\ge 1,$ is bounded and uniformly continuous, the
following estimate holds
\begin{equation}\label{2der}
|h(t,u)|\le n(r)|u|,
\end{equation}
where continuous function $n$ satisfies the conditions $n(0)=0,$ $n(r)\to 0,$ as $r\to 0.$

Let us introduce a function
\begin{equation}\label{Lyap}
S(t,u)=\left (\int\limits_{t}^\infty \varphi^2(s,t)ds\right )u^2 = s^2(t)u^2,\;\varphi(t,s)= \exp[\int_s^t a(\tau)d\tau].
\end{equation}
The improper integral converges uniformly w.r.t. $t$ and defines a bounded positive function $s(t)$ which is also
bounded away from zero. Indeed, one has an estimate
$$
\int\limits_{t}^\infty \varphi^2(s,t)ds\le C^2\int\limits_{t}^\infty \exp[-2\lambda(s-t)]ds = C^2/2\lambda.
$$
Denote $a_0 = \mathop{{sup}}\limits_t|a(t)|>0$. Then the estimate holds
$$
\int\limits_{t}^\infty \varphi^2(s,t)ds\ge \int\limits_{t}^\infty \exp[-2a_0(s-t)]ds = 1/2a_0.
$$
Thus for function $S$ the two-sided estimate is valid
\begin{equation}\label{tsided}
(1/2a_0)u^2 \le S(t,u)\le (C^2/2\lambda)u^2.
\end{equation}

Now let us show $S$ be a Lyapunov function for the equation (\ref{linear}) in a neighborhood of its solution $u=0$.
This means that for any initial point $(t_0,u_0)$
the following estimate for the derivative of $S$ along a solution $u(t)$ through $(t_0,u_0)$
to the equation (\ref{linear}) takes place
$$
\frac{d}{dt}S(t,u(t))|_{t_0}= \frac{d}{dt}[s^2(t)u^2(t)]|_{t_0}= 2s(t_0)s'(t_0)u_0^2 +
2s^2(t_0)u_0[a(t_0)u_0 + h(t_0,u_0)].
$$
Differentiation of $s(t)$ gives $s'(t)= -1/2s(t) - a(t)s(t).$ Hence we get
$$
2s(t_0)s'(t_0)u_0^2 + 2s^2(t_0)u_0[a(t_0)u_0 + h(t_0,u_0)] = - u_0^2 + 2s^2(t_0)u_0 h(t_0,x_0).
$$
Therefore one has
$$
\frac{d}{dt}S(t,u(t))|_{t_0}= -[1-2s^2(t_0)h(t_0,u_0)/u_0]u_0^2.
$$
Due to the estimate (\ref{2der}), the function within square brackets can be made less than $1/2$
choosing $|u_0|$ small enough.
Let us observe that the derivative of $S$ along the system at the point $(t_0,u_0)$ is the inner product $(S_t,S_u)$ of
$\nabla S$ and the tangent vector $(1,u')$ to the integral curve at the
point. From the estimate obtained it follows that at a fixed $c^2 > 0$ along the curve $S(t,u)=c^2$
the inner product in bounded away from zero uniformly in $t$. Fixing $c > 0$ small enough
gives an uniform neighborhood of the solution $u=0$ defined as $S < c^2.$

Now we consider again the equation (\ref{1d}). Henceforth we denote $S^1_0$ the section $t=0$ in $\R\times S^1.$
Suppose an integral curve $\gamma$ possesses exponential dichotomy
of the stable type on $\R_+$. Then all ICs through initial points being close enough to the point $\gamma(t_0)$
at $t=t_0>0$ possess also the exponential dichotomy of the same type. This implies a set of traces on
the section $t=0$ for all ICs with the exponential dichotomy of stable type on $\R_+$ be open
and therefore to consist of a collection of disjoint intervals (the cardinality of this collection of intervals is
maximum countable).
Thus, the set of points on the section $S^1_0$ through which ICs with exponential dichotomy
of the unstable type on $\R_+$ pass, form a closed set. The cardinality of this set is an interesting question
when the Assumption 1 holds.
Our further goal is to distinguish a class of NVFs that are rough and can be classified somehow.
This requires of more rigid conditions on the behavior of their ICs.  We shall call the set of ICs
with exponential dichotomy of the stable type on $\R_+$ by a {\em stable bunch} if all of
them have the same asymptotic behavior, i.e. any two solutions $x_1(t),$ $x_2(t)$ from this
stable bunch satisfy inequality $|x_1(t)-x_2(t)|\to 0 \;\mbox{\rm (mod\;1)}$, as $t\to
\infty$. In fact the union of all IC of the same stable bunch are a global
stable manifold for any IC from this bunch. The boundaries of a stable
bunch consists of one or two ICs which possess exponential dichotomy of
the unstable type on $\R_+.$

Similar properties are valid for ICs which possess the exponential dichotomy on semi-axis
$\R_-.$ Here we define an {\em unstable bunch} as the union of all
ICs that possess exponential dichotomy of the unstable type on semi-axis $\R_-.$
Then our next assumption is the following
\begin{assumption}.
There are finitely many stable bunches and finitely many unstable bunches.
\end{assumption}

The traces on $S^1_0$ of those ICs which belong to one stable bunch is an interval. Extreme points of
this interval correspond to ICs that possess exponential dichotomy of the unstable type of semi-axis $\R_+$.
Similarly, the traces of ICs from one unstable bunch on $S^1_0$ is an interval and its extreme points correspond to
those ICs which possess exponential dichotomy of the stable type on semi-axis $\R_-$.
Thus, from this Assumption it follows the existence of finitely many ICs with exponential dichotomy of the unstable
type of semi-axis $\R_+$ and also does for ICs with exponential dichotomy of the stable type on semi-axis $\R_-.$

Our third assumption is
\begin{assumption}. No solutions exist such that this solution possesses simultaneously exponential dichotomy
of the unstable type on $\R_+$ and exponential dichotomy of the stable type on $\R_-.$
\end{assumption}

One more property of the class of ODEs under study has to be discussed.
Consider some stable bunch and let $\gamma_1, \gamma_2$ be its boundary
solutions (left and right, in accordance with the orientation of $S^1$,
respectively). Hence these two solutions possess exponential dichotomy of
the unstable type on $\R_+.$ Choose some solution $\gamma_0$ of the bunch. Due
to the Hadamard-Perron theorem (or merely to the existence of a Lyapunov
function) one can choose an uniform neighborhood $U_0$ of $\gamma_0$ on $\R_+$
within the bunch whose boundary curves are uniformly transversal to solutions of ODE
through these boundary curves and these solutions enter to the
neighborhood and stay there forever. Similarly, for the solution
$\gamma_1$ there is its uniform neighborhood $U_1$ with the property of
uniform transversality of its boundary curves but all solutions through
this neighborhood leave $U_1$ at some $t$ specific for such solution.
Without loss of generality we suppose $U_0, U_1$ not intersecting and
we can choose them as thin as we wish. Now an important question arises:
consider for $t\ge 0$ the semi-strip between the right boundary curve of $U_1$ and left
boundary curve of $U_0$. Solutions within this semi-strip enter to it from the left and
leave it from the right boundary curve. These solutions cannot stay in this semi-strip for all $t>0$,
since all of them belong to one stable bunch.
Is the passing time through this semi-strip is
uniformly bounded for all such solutions? The answer to this question is
essential for constructing equimorphism which conjugates equivalent ODEs.
We shall present below an example of a nonautonomous ODE that shows the answer to
this question is generally negative. Therefore this property has to be
imposed to avoid a structural instability. To this end, let us choose a
sufficiently thin disjoint neighborhoods $\mathcal U_1,\ldots,\mathcal U_s$ of all ICs which possess exponential
dichotomy of the unstable type on $\R_+$. Also we select by one IC from
every stable bunch and choose their thin disjoint neighborhoods $\mathcal V_1,\ldots,\mathcal
V_r$. These all neighborhoods $\mathcal U_i,\mathcal V_j$ can be chosen
disjoint.
\begin{assumption}
For any $\eps > 0$ there are neighborhoods $\mathcal U_1,\ldots,\mathcal U_s, \mathcal V_1,\ldots,\mathcal
V_r$ such that the passage time of ICs from one boundary curve to another
one is bounded from above.
\end{assumption}
Now we shall call scalar ODEs on $S^1$ which obey Assumptions 1-4 being {\em gradient-like}.

To demonstrate the essentiality of the assumptions imposed,
we present examples of differential equations (\ref{1d}) showing that violating any of these three assumptions
may give an ODE being nonrough. Our first example was found in \cite{Ler_dis} and
shows that if there is a solution that fails to possess the exponential
dichotomy on some of two semi-axes $\R_+, \R_-$, then the equation may become nonrough and its structure
changes at a perturbation.
Consider a nonautonomous differential equation $v$ on the circle with
coordinate $\varphi \in [0,2\pi)\;\mbox{\rm\;mod}\;2\pi$
$$
\dot \varphi = \cos\varphi(\sin^2\varphi+e^{-t^2}).
$$
The equation is reversible w.r.t. the involution $L_1:\varphi \to 2\pi -
\varphi$ and the involution $L_2: \varphi \to \varphi + \pi.$ Recall, this means that if $\varphi(t)$
is a solution to the equation, then so does $L_i\varphi(-t)$, that is
$\varphi_1(t) = 2\pi-\varphi(-t)$ and $\varphi_2(t)=\varphi(-t)+\pi$.

As $t\to \pm\infty$ this equation tends to the autonomous differential equation $v_*$ given as
$\dot \varphi = \cos\varphi\sin^2\varphi$ that has a nonhyperbolic equilibrium $\varphi =0.$
The equation $v_*$ considered as nonautonomous has the following foliation onto
ICs. There are four stationary solutions $\varphi(t) =
0,\pi/2,\pi,$ $3\pi/2.$ Two of them, $\varphi =\pi/2,3\pi/2$ possess exponential
dichotomy on $\R$, the first one is asymptotically stable and the second
is asymptotically unstable. More two stationary solutions do not possess
 the exponential dichotomy neither on $R_+$ no on $\R_-.$
ICs corresponding to stationary solutions divide $S^1\times \R$ into four strips. All ICs
in each strip go from one stationary solutions (at $-\infty$) to another one (at $\infty$).

Now we take into account in the equation $v$ the rapidly decaying term $\exp(-t^2)$.
Solutions $\varphi = \pi/2,3\pi/2$ stay unchanged and have the same
type of dichotomies. Let us consider the invariant strip $-\pi/2\le \varphi \le
 \pi/2$ with its foliation into ICs. Note that $-\pi/2 = 3\pi/2\;(\mbox{mod\,}2\pi).$
One can prove
\begin{lemma}
There are two special ICs in this strip. One of them $\gamma_1$ is monotonically
increasing and tends to $\varphi = 0$ as $t\to \infty$ but tends to
$\varphi = -\pi/2$ as $t\to -\infty$. Another one $\gamma_2$
is also monotonically increasing and is defined by reversibility as
$\varphi_1(t)= -\varphi(-t)$ where $\varphi(t)$ is the solution for $\gamma_1$.
\end{lemma}

Two special solutions indicated in Lemma do not possess the exponential dichotomy on either $\R_+$ (for $\gamma_1$) or
on $\R_-$ (for $\gamma_2$). Observe that the straight-line $\varphi =0$ is
transversal to ICs intersecting it. All ICs between two special ones $\gamma_1, \gamma_2$ intersect
this straight-line and tend to
$\varphi =\pi/2$ as $t\to \infty$ and to $\varphi =-\pi/2$ as $t\to -\infty$
and thus possess the exponential dichotomy. Other ICs lying in the strip
below $\gamma_1$ or above $\gamma_2$ do not possess exponential dichotomy
on the related semi-axes since they tend in some direction in time to
nonhyperbolic solution ($\gamma_1$ or $\gamma_2$). The behavior of ICs in
the strip $\pi/2 \le \varphi \le 3\pi/2$ is similar.

\begin{figure}[h]
\begin{center}
\includegraphics[width=12cm, keepaspectratio,
scale=1.5]{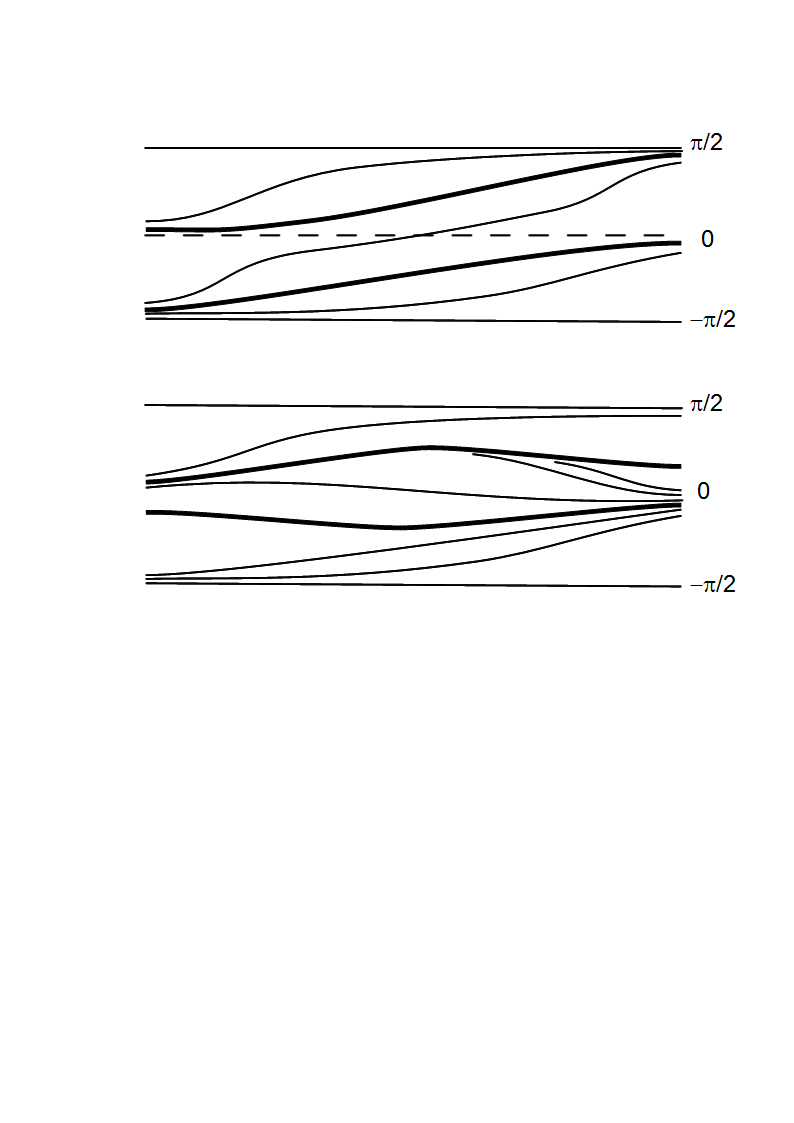}
\caption{Bifurcation from infinities, $\mu=0,$ $\mu > 0.$}
\label{1}
\end{center}
\end{figure}

This ODE under consideration changes its uniform structure under perturbation.
Indeed, consider a perturbed equation $\dot \varphi = \cos\varphi(\sin^2\varphi+\exp[-t^2]-\mu)$ with small positive $\mu$.
One can verify that two new bunches have appeared for small $\mu$, one stable and one unstable.
They arise because of the bifurcation, when nonhyperbolic equilibria on $\pm\infty$ bifurcate into
two close hyperbolic equilibria. The related foliations are shown in
Fig.\ref{1}. The structure and bifurcation presented is similar to that
observed in the Ricatti equation $\dot x = x^2 + \exp[-t^2]- \mu.$

\begin{figure}[!ht]
\begin{minipage}[h]{0.3\linewidth}
\center{\includegraphics[width=1\linewidth]{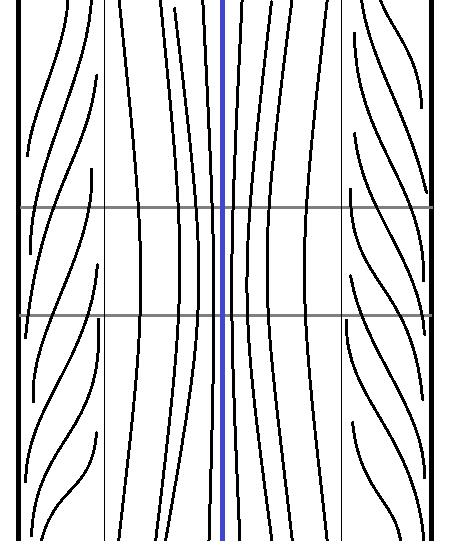} \\ a) Degenerate IC}
\end{minipage}
\hfill
\begin{minipage}[h]{0.5\linewidth}
\center{\includegraphics[width=1\linewidth]{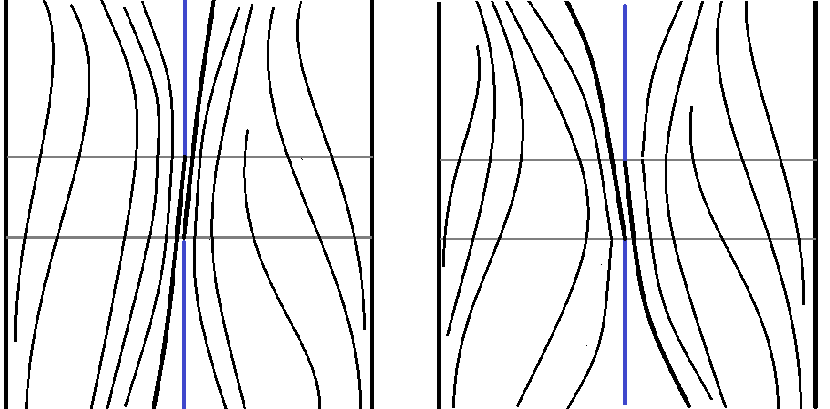} \\ b) Two types of bifurcation}
\end{minipage}
\caption{{\footnotesize Bifurcation at the violation Assumption 3.}}
\label{2}
\end{figure}

The second example shows that the violation of the third assumption leads
to a structurally unstable ODE. Consider an ODE that has a solution with exponential
dichotomy of the unstable type on $\R_+$ and let a point $m$ be its trace on the section
$S^1_0$. Suppose the ODE under consideration is such that IC
through $m$ possesses exponential dichotomy of the stable
type on $\R_-$. This implies that there is a sufficiently thin
uniform neighborhood of this IC on $\R$ within which the only IC
through $m$ lies in $U$ entirely for all $t$. If for this ODE two first
Assumptions hold, then point $m$ serves the boundary point of two
intervals corresponding two neighboring stable bunches. Their other
boundaries (to the left and to the right from $m$) on $S^1_0$
are points $m^+_1, m^+_2$ (these two points may coincide
if the closure of two intervals make up $S^1_0$).
Let choose by one IC from both stable bunches. Denote them as $\gamma_1,
\gamma_2.$ The minimal distance between related ICs on $t\ge 0$ is
some positive number $\rho.$ (Fig.\ref{2})

Now we perturb the ODE in such a way that the perturbed ODE would have one
IC with the dichotomy of the unstable type on $\R_+$ close enough on $\R_+$ to
former one through $m$ and another IC with the dichotomy of the stable type
on $\R_-$ being close enough on $\R_-$ to
former IC through $m$. This perturbation can be chosen in two ways as shown in
Fig.\ref{2}. Since the uniform equivalence preserve asymptotic behavior of
ICs, these two ODEs are not uniformly equivalent.

To present an example that demonstrates a possibility to be unbounded for the time of passage
from one boundary curve to another one in a transitory strip, we construct an
example of equation $\dot x = f(t,x)$ on the segment $I=[-1,1]$ instead of $S^1$.
It is easy to glue two such strips to obtain an example on $S^1.$
This shows the Assumption 4 be essential.

\begin{figure}
\centering
\parindent=0pt
\includegraphics[width=0.5\textwidth]{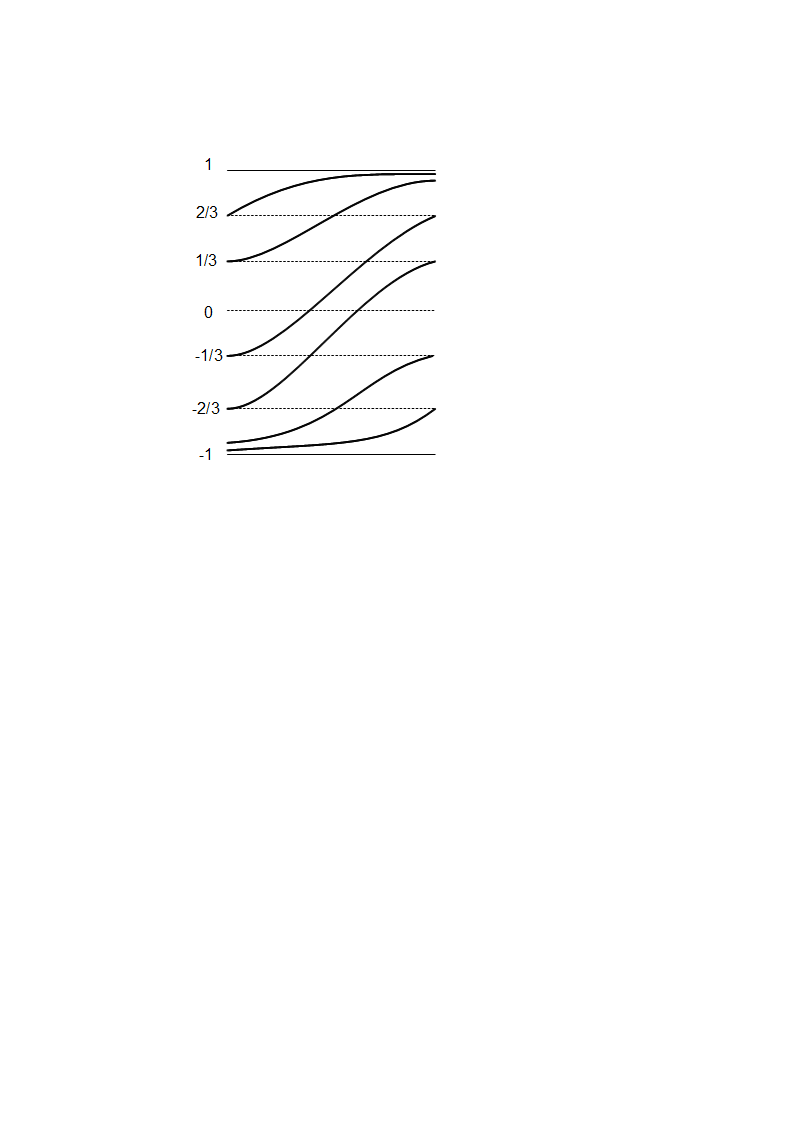}
\caption{Foliation inside the rectangle.}
\label{3}
\end{figure}

Consider the semi-bounded strip $\R_+\times I$ with coordinates $(t,x)$ and
choose an increasing sequence of numbers $t_n \to \infty$ such that $t_0 =0,$ $t_{2k}-t_{2k-1} =
1,$ $k\ge 1,$ and $t_{2k+1}-t_{2k} \to \infty,$ as $k\to \infty.$ Within
segments in $t$ given as $[t_{2k},t_{2k+1}]$, $k\ge 0,$ we take the function $f(t,x)= g(x)$
independent on $t$ and the same for all such segments. Function $g$ is
of class $C^\infty$, odd and defined as follows
$$
g(x)=\left\{\begin{array}{l}x+1,\;\mbox{if}\, x\in [-1,-2/3],\\g(x)> 0 \,{\rm positive\,with \,the \,only\,maximum\,as}\,
x\in (-2/3,-1/3),\\g\equiv 0,
\,\mbox{if}\, x\in [-1/3,0],\\g(x)=-g(-x)\; \mbox{for\, positive}\;x.\end{array}\right.
$$
We choose the function $f$ independent in $t$
(autonomous) for all $t\in \R_+$ in semi-strips $x\in [-1,-2/3]$ and $x\in [2/3,1]$.
Finally, in rectangles $[t_{2k-1},t_{2k}]\times I$, $k\ge 1,$ of the constant length
we defined $f(t,x)$ in such a way that $f$ is positive for $t\in (t_{2k-1},t_{2k})$ with the related foliation
into integral curves within such the rectangle as in Fig.\ref{3}.

Its specifics is that the solution $x= 1/3$ from the left region goes up and
intersects the line $x = 2/3$ between $t=t_{2k-1}$ and $t=t_{2k}$ but
solution $x = -1/3$ goes up and cut the line $x=2/3$ at $t= t_{2k}.$
Similarly, solutions $x= 1/3$ and $x=-1/3$ from the right region, as $t$
decreases, intersect the line $x=-2/3$ respectively at $t= t_{2k-1}$ and
between $t= t_{2k-1}$ and $t= t_{2k}$. All this guarantees that
\begin{itemize}
\item each
integral curve, except $x= -1$, intersects the line $x= 2/3$ entering thus to the strip
where the equation is $\dot x = - (x-1),$ such solution possesses
exponential dichotomy of the stable type on $\R_+$;
\item there are countably many solutions of the full equation which coincide with $x= -1/3$ on the
segment $[t_{2k},t_{2k+1}]$ have their passage time from $t=-2/3$ to
$t=2/3$ greater than $t_{2k}-t_{2k+1}$ and hence these times unbounded.
\end{itemize}

\section{Examples of gradient-like scalar ODEs}

Here we present several simple examples of scalar equations of the gradient
type on $S^1$. The first example is trivial, it is an autonomous scalar ODE $\dot x =
f(x)$ with a smooth $2\pi$-periodic $f$ having only simple zeroes ($f'(m)\ne 0$ at such a zero
$m$) and hence there are only finitely many zeroes (equilibria). Due to compactness of $S^1$ all solutions to ODE
tend, as $t\to \pm\infty$, to some of its equilibrium, their simplicity guarantees the exponential dichotomy
of the related constant solutions on $\R$, as well as the exponential dichotomy on the semi-axis $\R_+$ or $\R_-$
for all solutions that tend to this constant solution as $t\to \infty$ or $t\to -\infty$, respectively.

The second example is less trivial but also well known. Consider a
nonautonomous doubly periodic equation $\dot \varphi = f(t,\varphi)$. The
function $f$ is $2\pi$-periodic in both variables and smooth (in fact, it is sufficient $f$ to be continuous
in $t,\varphi$ and differentiable in $\varphi$ with a continuous derivative in $(t,\varphi)$). Then
the Poincar\'e map $P: S^1 \to S^1$ in $t$ is well defined. This map is a diffeomorphism (at least
$C^1$) and has a Poincar\'e rotation number. Suppose this number is rational $p/q$ with incommensurate $p,q$. Then
by the Poincar\'e theorem (see, for instance, \cite{Hartman}) this map has a periodic orbit of the period $q$
and any other orbit is either periodic one with the same period or tend
to a periodic orbit as $n\to \infty$ ($n$ is the number of iterations) and to a
periodic orbit (possibly to another one) as $n\to -\infty$. Suppose all these periodic orbits are
hyperbolic, i.e. $DP^q \ne \pm 1$ at any periodic orbit. Then there are only
finitely many periodic orbits and hence this nonautonomous ODE is of the gradient type.

The third example is a so-called asymptotically autonomous scalar ODE
\cite{Markus,SY,rassmus}. Suppose a differential equation on $S^1$ be given $\dot \varphi =
f(t,x)$ such that $f \to f_+$ as $t\to\infty$ and $f \to f_-$ as $t\to -\infty$
where the limit is taken in the space $V^r(S^1),$ $r\ge 1.$ We assume
that both ODEs on the circle with functions $f_+, f_-$ have only simple
zeroes as in the first example. Then any IC of the given equation asymptotically
approaches as $t\to \infty$ to some $\varphi_k$, $f_+(\varphi_k)=0$, and asymptotically
approaches as $t\to -\infty$ to some $\varphi_j$, $f_-(\varphi_j)=0$ \cite{Markus}. It
follows from here that any integral curve possesses an exponential
dichotomy both on $R_+, R_-.$ In order the equation would be gradient-like,
one needs to require in addition that no IC exists that asymptotically tends as $t\to -\infty$ to a stable zero
of $f_-$ and as $t\to \infty$ does to an unstable zero of $f_+$.
Then this ODE is a nonautonomous  gradient-like one. As
an example of such the equation one can construct a glued differential equation. We take a
rough scalar autonomous equation with the function $f_+$ and another rough
scalar autonomous equation with the function $f_-$. The nonautonomous
equation is glued as follows. Let for $t\le -T$ the equation with $f_-$ is
considered, and for $t\ge T$ does the equation with $f_+.$ In the layer between two
sections $-T\le t\le T$ we can specify a smooth foliation in such a way that traces
of ICs for the equation with $f_-$ would join with traces of ICs to the equation with $f_+$
and all curves of the foliation in this layer would intersect transversely sections $S^1_t$,
$|t|\le T$. This makes the complete foliation smooth.
The derivatives in $t$ for ICs will give gradient-like nonautomous ODE.

In the similar way examples of a nonautonomous asymptotically periodic
scalar gradient-like ODEs can be constructed. To that end, one needs
to choose two periodic limiting ODEs as in the example 2 above (they play the role of
limiting equations at $\pm\infty$) given by functions $f_-(t,x)$ and $f_+(t,x)$
being periodic in $t$ with possibly different periods. Suppose their
related foliations into ICs are disposed lower of $S^1_{-T}$ (for $f_-$) and above of
$S^1_{T}$ (for $f_+$). Then we get on $S^1$ two 1-periodic in $x$ functions $f_-(-T,x)$ and
$f_+(T,x)$ being two points in the Banach space $V^1(S^1).$ Let us join
these points by a smooth compact path with the parametrization $t\in [-T,T]$.
Then we get a glued function $f(t,x)$ which defines the nonautonomous
ODE on $S^1$ and coincides with $f_-(t,x)$ when $t<-T$ and does with $f_+(t,x)$ when $t>T.$
The only thing we need to care about is that in a middle part of $\R$
traces of stable periodic solutions to negative limiting ODE (for $t< -T$)
would not be connected with periodic solutions of positive limiting ODE being unstable
hyperbolic for $t>T$. This is achieved by a small perturbation of the ODE in the
layer.

\section{The properties of gradient-like NVFs on $S^1$ and u-invariant}

In this section we define an invariant which is able to recognize
different gradient-like ODEs. This invariant is of combinatorial type,
i.e. it is defined by the finite set of ingredients.

To construct the invariant we begin with some assertions necessary for the
construction.
\begin{lemma}\label{ex}
For a gradient-like scalar ODE at least one solution with exponential dichotomy of the
unstable type on $\R_+$ always exists, as well as at least one solution with exponential
dichotomy of the stable type on $\R_-$.
\end{lemma}
\proof Suppose, to the contrary, a gradient-like ODE has not any solution
with exponential dichotomy of the unstable type on $\R_+.$ This implies
that there exists the unique stable bunch and the union of traces on $t=0$
of ICs from this bunch is the whole $S^1_0.$ Choose some IC $\gamma$ from this
bunch. There is an uniform neighborhood $U$ of $\gamma$ such that all ICs
starting within $U$ stay inside it for all $t\ge 0.$ Thus, one can choose
such a neighborhood for any IC of the bunch. The union of these neighborhoods
covers $\R_+\times S^1$, thus the intersection of these neighborhoods with the
section $S^1_0$ gives a cover of $S^1_0.$ Compactness of $S^1_0$ allows
one to select a finite cover and consequently a finite set of ICs
$\gamma_1,\ldots,$ $\gamma_n$ around which these neighborhoods have been
chosen. Denote related traces on $S^1_0$ of these ICs as points $p_1,\ldots, p_n$
enumerated in accordance with the orientation of $S^1_0$. Take some positive number $\eps < 1/2n.$
For any $p_i$ there is $T_i > 0$ such that for $t>
T_i$ the distance between traces of those ICs, which start at $t=0$ at extreme
points of $U\cap S^1_0$, on the section $S^1_t$ will
be lesser than $\eps.$ Thus for $t > \mbox{\rm max\,}\{T_1,\ldots T_n\}$
the distances between traces of any two ICs from the bunch will be lesser
than $n\cdot \eps < 1/2.$ On the other hand, they should cover $S^1_t.$
This contradiction proves the lemma. The similar proof holds true for the
dichotomy of the stable type on $R_-.$
$\blacksquare$

Now we shall define a combinatorial invariant which distinguishes two uniformly
nonequivalent gradient-like ODEs on $S^1$. Let a gradient-like ODE be given.
Take the section $S^1_0$. We assume that $S^1$ is oriented by the coordinate. On the section
we consider the traces of all ICs which have
the dichotomy of unstable type on $\R_+.$ Denote this set of points as $u_1,\ldots, u_n$
where the order from $1$ to $n$ is defined by the orientation of $S^1_0.$
Similarly, we consider the traces of all ICs which have
the dichotomy of stable type on $\R_-.$ Denote them as $s_1,\ldots, s_m$
where the order from $1$ to $m$ is also defined by the orientation of $S^1.$
These two sets of points on $S^1_0$ do not intersect in accordance with
the Assumptions 1-3. Due to Lemma \ref{ex}, integers $n,m$ are both positive.
We shall call the set of points $\{u_1,\ldots,u_n,s_1,\ldots,s_m\}$ the
{\em equipped set of points}.
\begin{definition}
Two equipped sets on $S^1$ will be called equivalent, if there is
a homeomorphism $h:S^1 \to S^1$ such that $h$ sends the set
$u_1,\ldots,u_n$ of the first equipped set onto the set $u'_1,\ldots,u'_n$ of the
second equipped set, and the set $s_1,\ldots,s_m$ of the first equipped set to
$s'_1,\ldots,s'_m$ of the second equipped set.
\end{definition}
A class of equivalent equipped sets will be called u-invariant.
The main result of the paper is the following theorem
\begin{theorem}\label{equi}
Two nonautonomous scalar gradient-like ODEs on $S^1$ are
uniformly equivalent, iff they have the same u-invariant.
\end{theorem}
It is clear that two uniformly equivalent ODEs on $S^1$ have the same
u-invariant. The proof of the inverse assertion is the main task.

Gradient-like ODEs possess also the property of structural stability.
This notion was introduced by Andronov and Pontryagin in 1937 \cite{AndP} for autonomous vector fields on
2-dimensional sphere (or on a disk with the proper assumptions on the orbit behavior
on the boundary cicrle) and was studied in many details, first on sphere
(Leontovich-Andronova, Maier \cite{LM}) and then on the torus (Maier). The final
result for any smooth closed 2-dimensional manifold was done by Peixoto
\cite{Peix}. Multidimensional case was studied first by Smale \cite{Sm}
who defined the class of vector fields and diffeomorphism which were
called later Morse-Smale systems. Their study lasts till now \cite{Grines,GMP}.

Leaning on Theorem \ref{equi} we will also prove
\begin{theorem}
A gradient-like NVF $v$ on $S^1$ is rough, that is, for any $\eps > 0$
there exists $\delta(\eps)$ such that all NVFs
in the $\eps$-neighborhood of $v$ in the space $\R\to V^1(S^1)$ are uniformly
equivalent to $v$ and an equimorphism realizing
this equivalence can be chosen $\delta$-close to the identity mapping $id_{\R\times S^1}$.
Moreover, all shifted NVFs $v_{t+\tau}$ are also rough with the same $\delta.$
\end{theorem}
Proofs of both theorems will be outlined in the next two sections.

\section{A method to construct an equimorphism}\label{method}

To demonstrate the method of constructing an equimorphism for two gradient-like ODEs with the same u-invariant,
we present here first the simplest possible situation
when both differential equations are linear. It gives, in a sense, the simplest uniform version
of the Grobman-Hartman theorem \cite{Hartman}. For the case of two
gradient-like nonlinear ODEs the local construction is similar.

Consider two scalar linear homogeneous differential equations
\begin{equation}\label{ls}
x'=a(t)x(t),\;y'=b(t)y(t).
\end{equation}
Denote $u(t,\tau) $ the solution of (\ref{ls}) with the condition $
u(\tau,\tau)=1,$ and let $u_{1}(t,\tau) $ be the similar solution for the
second equation
$$u(t,\tau)=\exp[\int\limits_{\tau}^{t} a(s)ds],\;
u_{1}(t,\tau)=\exp[\int\limits_{\tau}^{t} b(s)ds].
$$
We assume both these linear ODEs to possess exponential dichotomy of the stable type on $\R_+$
\begin{equation}
u(t,\tau)= \exp[\int\limits_{\tau}^{t}a(s)ds]\le M \exp[-\lambda (t-\tau)]; \ t\ge\tau,\ \lambda >0,\ M > 0,
\end{equation}
and the same holds for the second equation with constants $M_1,\lambda_1.$

We defined above a Lyapunov function (\ref{Lyap}) for such equation. For
equations under consideration we denote related functions as $S(t,x), s(t), u(t,\tau),$
$S_1t,x), s_1(t), u_1(t,\tau).$

From expressions for $u,u_1$ the inequalities follow
\begin{equation}
\dfrac{1}{2a_{0}} \le {s^2(t)} \le \dfrac{M^2}{2\lambda}
\end{equation}
\begin{equation}
\dfrac{1}{2b_{0}} \le {s_{1}^2(t)} \le \dfrac{M^2_{1}}{2\lambda_{1}}
\end{equation}
where $a_{0}= \sup\limits_{t}|a(t)|,$ $b_{0}=\sup\limits_{t} |b(t)|.$

Observe that functions $s(t),s_{1}(t)$ obey the equations
$$\dfrac{ds}{dt}=\dfrac{-2a(t)s^2(t)-1}{2s(t)}$$
$$\dfrac{ds_{1}}{dt}=\dfrac{-2b(t)s_{1}^2(t)-1}{2s_{1}(t)}.$$

Our goal in this section is to prove that the map from the semi-strip $D$
in the plane $(x,t)$ around $x=0$ onto the semi-strip
$D_1$ in the plane $ (y,t)$ around $y=0$ is an equimorphism. The strip $D$ is defined by its boundary curves $x=0$
and $x = C^{*}/s(t)$, $t\ge 0$, where $C^*$ is some positive constant. This semi-strip
contains the right curve but does not contain the left one $x=0$. Similar formulae define
the strip $D_1$ with the change $x\to y,$ $C^* \to C_1^*,$
$s(t) \to s_1(t),$ $t\ge 0$.

Let us define the map $\Phi: D \to D_1$.
First we change the coordinates in both semi-strips. Instead of $(x,t)$ we take as new coordinates $(C,t)$,
similarly, we change $(y,t)$ to $(C_1,t).$
In accordance with the formulae for Lyapunov functions, these coordinate transformations are smooth in both variables
and linear in $x$ (respectively, $y$) with uniformly bounded coefficients being uniformly bounded away from zero.
Then the map $\Phi$ is given as follows. Take any point $(C,\tau)\in D$ and consider the
integral curve of the first equation through this point. This curve
intersects transversely the level line of the Lyapunov function $S$ defined by $C^*$ at some time
$T(C,\tau) \le \tau$. We fix some smooth monotone function $\alpha: \mathbb R \to \mathbb R$
given as $T_1 = \alpha(T)$ with $\alpha'$
positive uniformly bounded from above and from zero: $0< r_0\le \alpha' \le
r_1.$ For instance, as a such map one can take the identical map $\alpha(T)= T.$ To facilitate
the exposition, we take just this identical map. Then the point on the right boundary of $D$ with coordinates
$(C^*,T(C,\tau))$ is transformed at the point $(C_1^*, T(C,\tau))$ on
the right boundary curve of $D_1.$ On the integral curve $(t,y(t))$ of the second equation through
this latter point we choose that its point whose $t$-coordinate is $\tau$. Then $C_1$-coordinate
of this point is $C_1 = y(\tau)s_1(\tau).$

Integral curve through the point $(C=s(\tau)x_0,\tau) $ is
$$x(t)= x_{0}\exp[\int\limits_{\tau}^{t}a(u)du] =  \dfrac{C}{s(\tau)}
\exp[\int\limits_{\tau}^{t}a(u)du],\;x_0 >0. $$
This curve intersects the boundary level line $s(T) x = C^{*}$ at the point with
coordinates $(C^{*}, T)$ where $T$ is found from the
equation
\begin{equation}
C s(T)\exp[\int\limits_{\tau}^{T}a(u)du] - C^{*} s(\tau) = 0.
\end{equation}
Denote $R(T,C,\tau)$ the left hand side of this equation. Function
$T(C,\tau)$ is a solution to the equation $R(T,C,\tau)=0$ existing by the implicit function theorem.

Integral curve $(t,y(t))$ of the second equation through the point $(C_{1}^{*},T)$
is given as (its second coordinate)
 $$ y(t) = \dfrac{C_{1}^{*} \exp[\int\limits_{T}^{t}b(u)du]}{s_{1}(T)}.$$
For the point $\Phi(C,\tau)$ the time variable on this curve is $\tau_{1}=\tau.$
$C_{1}$-coordinate of this point is found from the equation
\begin{equation}\label{c1}
 C_{1}= \dfrac{C_{1}^{*} s_{1}(\tau)\exp[\int\limits_{T}^{\tau}b(u)du]}{s_{1}(T)}.
 \end{equation}

The right hand side of the equation (\ref{c1}) will be denoted as
$g(T,\tau).$ Thus the map $\Phi$ is given in coordinates as follows
\begin{equation}\label{pxi}
C_1 = g(T(C,\tau),\tau),\; \tau_1 = \tau.
\end{equation}
We need to prove the uniform continuity of this map along with its inverse one.

In fact, we need only to prove the uniform continuity of $g(T(C,\tau),\tau)$ in the strip $0<C\le C^*. $
The uniform continuity of the inverse map is proved in the same way, since
the inverse map is defined in the similar way. The uniform continuity of
$\Phi$, when it has been proved, allows to extend the map uniquely to the closed strip
$0\le C \le C^*$ by the relation $\tau_1 = \tau.$

The uniform continuity of $g$ is proved in two steps because the derivatives of $\Phi$ is not
bounded if we consider the full strip $0< C\le C^*$. Therefore we fix some $v > 0$
and prove that $\Phi$-pre-image in $D$ of the line $C_1 = v$ is a curve being the graph of a function
$\varphi(\tau),$
and this function obeys estimates $0<d_1\le \varphi(\tau)\le d_2 < C^*.$
Moreover, if $v\to +0$ then both $d_i(v)$ tend also to zero. Thus choosing
$v$ small enough we get that the semi-strip $0< C \le d_1$ is transformed
on the semi-strip in $D_1$ which belongs to the semi-strip $0<C_1\le v.$
Since $\tau$-coordinate preserves by $\Phi$ this gives uniform continuity
near the line $x=0$. In the closed semi-strip $d_1\le C \le C^*$ the
derivatives of $\Phi$ are uniformly bounded and this guarantees the
uniform continuity for $\Phi$ there.

Now consider in $D$ pre-images of the Lyapunov function $S_1$ level lines defines by constants $v$
in $D_1$. First we will show that this curve is given as $C =
\varphi(\tau)$ for which the inequalities hold
$$
0\le d_1 \le \varphi_2 \le d_2.
$$
and $d_1, d_2$ tend to zero, if $v \to +0$.

To this end, we solve the equations (\ref{pxi}) w.r.t. $C$ setting $C_1 =
v$. Since the value of $\tau$ preserves, then the
solutions will give a function $\varphi(\tau)$ and we need only to
estimate them. These functions are given in the following way
$$
\varphi_i(\tau)= \frac{C^*s(\tau)\exp[\int\limits_{T}^\tau a(u)du]}{s(T)},\; T = T(v,\tau).
$$
Here function $T(v,\tau)$ is a solution of the equations (\ref{c1})
with $C_1 = v$. From this equation an estimate for the
difference $\tau - T$ is derived using the exponential dichotomy
$$
\tau - T \le \frac{1}{\lambda_1}\ln\frac{M_1C_1^*s_1(\tau)}{v s_1(T)}\le
\frac{1}{\lambda_1}\ln\frac{M_1C_1^*\sup(s_1)}{v \inf(s_1)} = A - \frac{1}{\lambda_1} \ln(v).
$$
We consider $v$ small enough positive, so r.h.s. will tend to infinity
as $v \to +0,$ but it is finite, if $v$ holds fixed. We have a similar
estimate from below
$$
\tau - T \ge A_0 - \frac{1}{b_0}\ln(v),\; A_0 = \frac{1}{b_0}\ln\frac{C_1^*\sup(s_1)}{\inf (s_1)}.
$$
Using these estimates for $\tau - T$ we come to the following estimates
for function $\varphi(\tau)$
\begin{equation}\label{var}
B_0(v)^{a_0/b_0}\le \varphi(\tau)\le B_1(v)^{\lambda/\lambda_1}.
\end{equation}
These estimates say that the pre-image of any strip $v\le C_1 \le C_1^*$ in
$D_1$ is a strip in $D$ lying in the the strip $0< d\le C \le C^*$. This
implies that for proving uniform continuity one may use estimates for the derivatives of $g$.
These derivatives are as follows
$$\dfrac {\partial g}{\partial C}= \dfrac{\partial g}{\partial T}\cdot\dfrac{\partial T}{\partial C},$$
and using (\ref{c1}) we get
\begin{equation}\label{der}
\dfrac {\partial g}{\partial T}= \dfrac{C_1^*s_1(\tau)\exp[\int\limits_T^\tau b(s)ds]}{s_1^3(T)}= \dfrac{C_1}{s_1^2(T)},\;
\dfrac {\partial T}{\partial C}=\dfrac{2s^2(T)}{C}.
\end{equation}
Therefore, the estimate for $\dfrac {\partial g}{\partial C}$ in the strip $C\ge d$
is given as
\begin{equation}\label{fin_der}
|\dfrac{\partial g}{\partial T}\cdot\dfrac{\partial T}{\partial C}|\le \dfrac{2C^*_1\sup(s^2)}{d\inf(s_1^2)} = R_1.
\end{equation}
Now the calculation of derivative
$$\frac{D g}{D \tau} = \frac{\partial g}{\partial T}\frac{\partial T}{\partial \tau} + \frac{\partial g}{\partial
\tau}$$
gives the formulae
$$
\frac{\partial g}{\partial T}= \dfrac{C_1}{s_1^2(T)},\;\frac{\partial g}{\partial \tau}= - \dfrac{C_1}{2s_1^2(\tau)},\;
\frac{\partial T}{\partial \tau}= \dfrac{C^*s^2(T)}{s^2(\tau)}.
$$
Thus, we come to the estimate
\begin{equation}\label{der_tau}
|\frac{D g}{D \tau}|\le |\frac{\partial g}{\partial T}||\frac{\partial T}{\partial \tau}|+|\frac{\partial g}{\partial
\tau}|\le C_1^*(C^*\dfrac{(\sup(s^2))}{\inf(s_1^2)\inf(s^2)}+ \dfrac{1}{2\inf(s_1^2)})= R_2.
\end{equation}
These estimates for derivatives show the uniform finiteness of derivatives.

Now we are able to complete the proof of the uniform continuity of $\Phi.$
Choose the metrics in both strips:
$$
\rho((C,\tau),(C',\tau'))= \max\{|C-C'|,|\tau - \tau'|\},\;
\rho_1((C_1,\tau_1),(C_1',\tau_1'))= \max\{|C_1-C_1'|,|\tau_1 - \tau_1'|\}.
$$
Let a positive $\eps$ be given. Put $v_1=\eps/2$ in (\ref{var}) and take
$\delta(\eps)= \min\{\eps/2,B_0(\eps/2)^{a_0/b_0}, \eps/2R_1, \eps/R_2\}$.
Take two arbitrary points $P=(C,\tau),$ $P'=(C',\tau')$ in $D$ at the distance
$\rho(P,P')< \delta.$ Three cases are possible: i) both $P,P'$ belong to the sub-strip $0 < C \le \delta;$ ii)
both points belong to the sub-strip $\delta \le C \le C^*;$ iii) one point,
say $P$, belongs to first sub-strip but $P'$ belongs to the second
sub-strip. In the first case both points $\Phi(P), \Phi(P')$ belong to the strip $0<C_1\le
\eps/2$. Since $\rho (P,P')= \max\{|C-C'|, |\tau - \tau'|\}\le \eps/2$ and $\Phi$ preserves $\tau$-coordinate,
then $\rho_1(\Phi(P),\Phi(P'))\le \eps.$ The same estimate holds for the
third case due to estimates for derivatives of $g$. The second case is
reduced to the combination of the first and third ones if we joint points
$P,P'$ by the segment on the right boundary of the semi-strip $0 < C \le \delta.$

The construction gives an equimorphism only in semi-strips $x>0$ and
$y>0.$ In order to get an equimorphism of the whole strip around $x=0$ and
$y=0,$ respectively, one needs to conform them along the straight-lines $x=0$ and
$y=0,$ respectively. In the construction presented it is done automatically since the
section $t=\tau$ is transform to the section $t=\tau.$ For the general
case one needs to care about it.

\section{Uniform equivalence and structural stability}

To prove Theorem \ref{equi}, we need to construct an equimorphism $h:\R\times
S^1 \to \R\times S^1$ which transforms the foliation into ICs of the first
gradient-like ODE to that of the second ODE.
The idea of the construction is close to that performed in the linear case (see
above) but with some modifications exploiting the structure of ODEs. The
main idea is to construct first an equimorphism in some neighborhood of
those ICs whose traces on $S^1_0$ enter to the equipped set (we call the set of these ICs a {\em skeleton} of the ODE)
and after that to extend this equimorphism to the whole manifold $\R\times
S^1.$ When an ODE under consideration is gradient-like one, then its
skeleton depends continuously on a perturbations in the class of uniformly
continuous bounded maps $\R \to V^r(S^1).$ This allows to prove that a
conjugating equimorphism can be found to be close to the identity map $\R\times
S^1.$ In the general case of two gradient-like ODEs with the same u-invariant the
construction of a conjugating equimorphism is geometrical and requires a
detailed description. We present here only some details of this
construction. A complete exposition will be done elsewhere.

We start with constructing a sufficiently thin neighborhood
of the skeleton of the ODE. As a skeleton we take the following collection of
ICs. Consider the equipped set of the ODE. Its $u$-points $(u_1,u_2,\ldots,u_n)$ on $S^1_0$
are distributed on $m$ intervals $(s_1,s_2),$ $(s_2,s_3),\ldots,$ $(s_m,s_1)$ in $S^1_0$
formed by the traces of $s$-points (Assumption 3). ICs through $u$-points which belong to the same
interval $(s_i,s_{i+1})$ are in the same unstable bunch. Choose IC $\gamma_1$ through the
point $u_1$. This IC belongs to some unstable bunch.
The intersection of this bunch with $S^1_0$ defines some interval whose boundary points are $s_i,$
$s_{i+1}$. This implies $\gamma_1$ to possess exponential dichotomy of
the unstable type on $\R$ (we call such IC to be globally unstable).
Recall that a sufficiently thin neighborhood of such IC does not contain wholly
any IC other than $\gamma_1$ and boundary curves of this neighborhood
are uniformly transversal to all ICs passing though points of these curves. Here there is
an alternative: either no other points from $u_2,\ldots,u_n$ exist which belong to
the same interval $(s_i,s_{i+1})$ or there are several such ICs. In the first case
the second ODE also has a corresponding IC $\gamma_1'$ which passes through the point $u_1'$ of its equipped set
and its uniform isolating neighborhood $U_1'$.
We can construct equimorphism from $U_1$ to $U_1'$ by the method similar to in the preceding section.

\begin{figure}[h]
\begin{center}
\includegraphics[width=15cm, keepaspectratio,
scale=1.0]{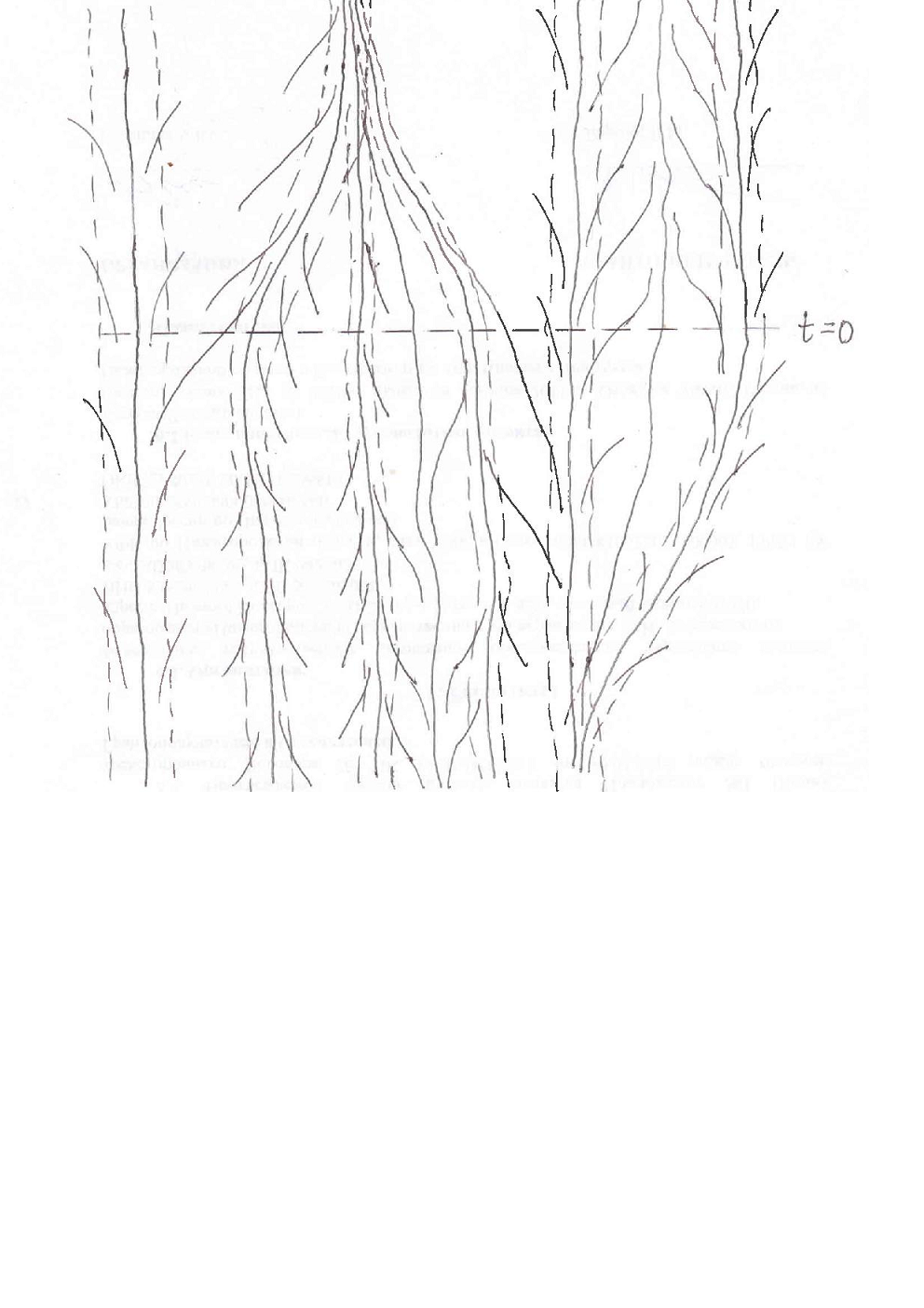}
\caption{Neighborhoods.}
\label{4}
\end{center}
\end{figure}

In the second case all related ICs
belong, as was said, to the same unstable bunch and they are asymptotic to each other as $t\to -\infty.$
Thus we have to construct a neighborhood $\mathcal U_1$ for all this collection of
ICs from the bunch. All these ICs are globally unstable ones. We choose some isolating
uniform neighborhood of IC through point $u_1$ and consider IC $\gamma_2$ through the point $u_2.$
There is its isolating neighborhood $U_2,$ in fact, many of them being as thin as we need.
IC $\gamma_2$ intersects one boundary curve of $U_1$ transversely at some point at $t < 0$. Thus we
can choose $U_2$ thin enough in order both boundary curves of $U_2$ would
intersect the same boundary curve of $U_1$ transversely and therefore
$U_2$ will belong to $U_1$ below some $t_1$ (see Fig.\ref{4}). Now we
proceed in the same way in order to construct a tree-shape neighborhood $\mathcal U_1$ for all
collection from one unstable bunch through $(s_i,s_{i+1})$ (see Fig.\ref{4}).

The boundaries of the neighborhood $\mathcal U_1$ are uniformly transverse to ICs but may contain a finitely many angular
points where two boundary curves from different $U_i$ adjoin. At such a point the IC through it is transverse to both local
smooth boundary curves. This IC possesses exponential dichotomy of the
stable type on $\R_+$. We shall distinguish different boundary curve of a
tree-shape neighborhood $\mathcal U_1$. First, there are two infinite in
both directions of time left $L$ and right $R$ boundary curves. We call them to be
outer ones. To specify all boundary curves we call the initial
neighborhood for $\gamma_1$ to be a trunk and other $U_i$ attached to the
trunk as branches. Then $L$ is the left boundary of the trunk drawing from $-\infty$
till its intersection with the left boundary curve of the leftmost branch (or till
$\infty$, if the trunk does not have branches to the left).
Analogously, we get the right boundary curve $R$ made up of the left
boundary curve of the trunk till its intersection with the right boundary
curve of the rightmost branch. Other boundary curves are of the $V$-shape
or of the semi-strip shape and formed by other boundary curves of inner
branches, they are semi-infinite and are extended till $\infty$.

Now we proceed in the similar way constructing neighborhoods $\mathcal
U_2,\ldots,$ $\mathcal U_s$ involving points
$u_{k_1+1},\ldots,u_{k_1+k_2},$ related with another unstable bunch, and so forth, till
we exhaust all collection of $u$-points in the equipped set. It is clear that the neighborhoods
can be constructed not intersecting and situated at a finite distance from each other.

\begin{figure}[h]
\begin{center}
\includegraphics[width=15cm, keepaspectratio,
scale=1.0]{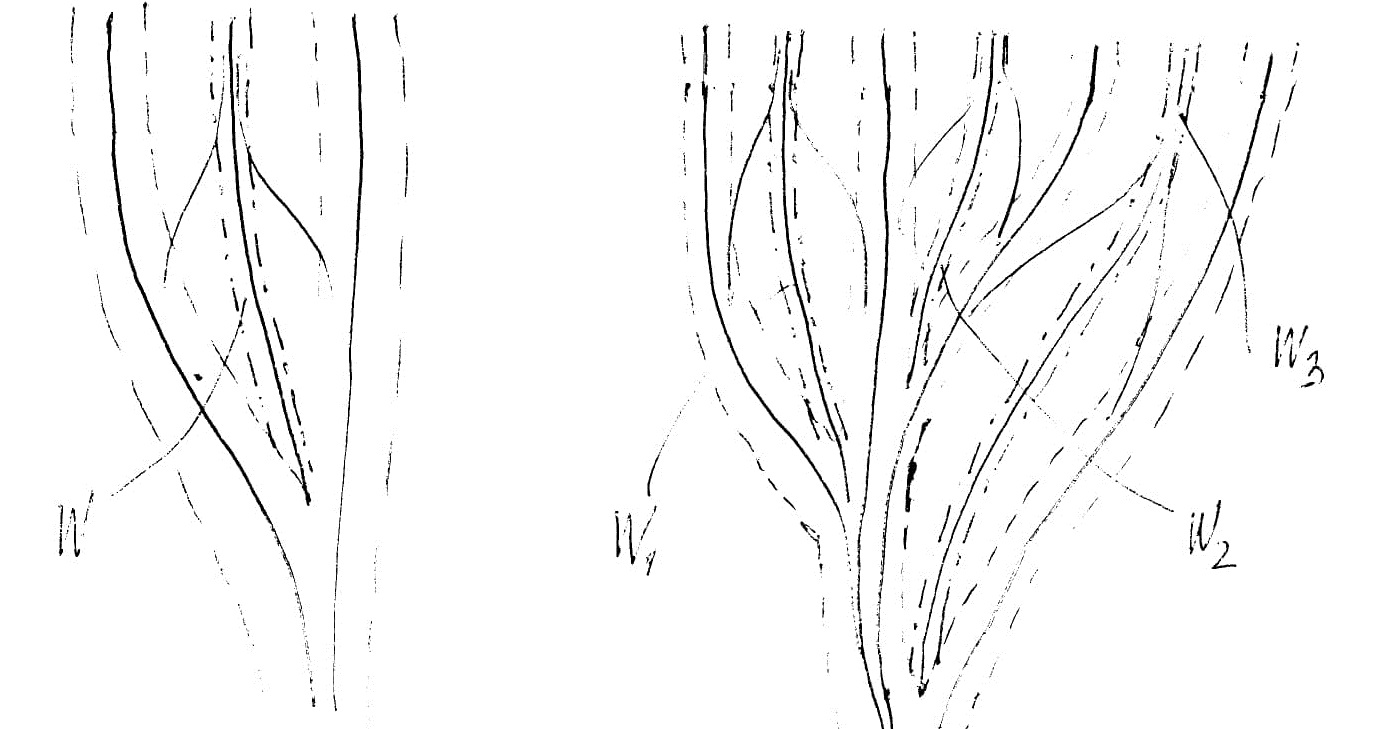}
\caption{Choice of W}
\label{5}
\end{center}
\end{figure}

Let us note that this construction does not involve so far $s$-points
$s_1,\ldots,s_m$. Now we continue the construction proceeding
in a similar way but starting with ICs through points $s_1,\ldots,s_m$ on
$S^1_0$ and using stable bunches with traces $(u_1,u_2),$ $(u_2,u_3),\ldots,$ $(u_n,u_1),$.
Thus, we construct the neighborhoods $\mathcal
V_1,\ldots,$ $\mathcal V_r$. These neighborhoods are also chosen sufficiently thin in order
to be separated from each other and from neighborhoods $\mathcal U_i.$ The $\mathcal U$-neighborhoods
constructed have the following characterization. They contain wholly only those
ICs which make up the skeleton. All other ICs either leave the related
neighborhood (for instance, for $\mathcal U_i$) when $t$ increasing at some finite
time. After that such an IC either in finite time enters to some $\mathcal
V_j$ through its lateral boundary curve (leftmost or rightmost one)
and this time is uniformly bounded from above for all such ICs (Assumption
4) or such an IC belongs to some stable bunch. In the second case we
choose in addition an IC in this stable bunch through angular point on the
boundary of $\mathcal U_i$ and select some its thin uniform neighborhood
$W^s_p$ whose boundary curves are transverse to the boundary curves of $\mathcal
U_i$ (see Fig.\ref{5})

For a neighborhood $\mathcal V_k$ we have a similar assertions with reversing
time. In particular, here one needs to choose in addition the
neighborhoods of ICs in unstable bunches which we denote as $W^u_l$. The
passage time from the boundaries of all neighborhoods constructed are
uniformly bounded due to Assumption 4.

\subsection{Conjugating equimorphism}

The construction of an equimorphism $h$ transforming one foliation into ICs to another one
is performed in two steps. At the first step we
define this equimorphism inside of the constructed neighborhoods $\mathcal
U_1,\ldots,$ $\mathcal U_s$, $\mathcal
V_1,\ldots,$ $\mathcal V_r$. At the second step we extend the
equimorphism outside of the union of these neighborhoods defining it wholly on $\R\times
S^1$. Since the construction is geometrical, it requires more room than we
have here. So, we present only some details, the complete exposition will
be done in the paper under preparation.

The basic construction of an equimorphism is performed in a thin semi-strip $t\ge T$ or $t\le -T$ around
fixed IC of the skeleton. Here $T$ is chosen large enough in order all such semi-strips do not
intersect each other. It can be also a thin neighborhood of an IC from a
stable (then $t\ge T$) or an unstable bunch (then $t\le T$). A
construction of an equimorphism within this semi-strip is the same as in
the previous Section \ref{method}. It is determined, if one sets an
equimorphism on the lateral boundary curve and a homeomorphism on the
segment of the bottom (or top) of the semi-strip. Here the levels of the
related Lyapunov function plays the main role. The lateral boundary curves of the semi-strips
themselves are the level lines of such a function. It requires a
neighborhood to be thin enough. We assume below that this construction is given.

It is clear that impossible to set equimorphisms on boundary curves
independently for different semi-strips. There are ICs which intersect
boundary curves of two different
semi-strips and this should be taken into account when constructing.
To verify easily that an equimorphism within a semi-strip is extended till
the basic IC of this semi-strip we shall set them identically in time for
$|t|\ge T.$ The construction the Section \ref{method} guarantees that
equimorphism within the semi-strip transforms the segment $t=c$ to the
segment $t=c$. In the cylinder $[-T,T]\times S^1$ this by means some
homeomorphism respecting the foliation.

In order to demonstrate the main points of the construction and how to
conform them to guarantee that any IC of the first ODE would transform by
the equimorphism to some IC of the second ODE we discuss two point.

The first point is constructing an equimorphism $h$ within a fixed
neighborhood $\mathcal U_i.$ To simplify notations, we omit the index $i$.
We start with the leftmost boundary curve $L$ of the neighborhood. To the right
from this curve the first IC from the skeleton is situated. We denote it
$\gamma_1$ (changing the numeration if necessary). We first construct the
equimorphism in the strip $\Pi_1$ between these two curves $L$ and $\gamma_1.$
Observe that on $L$ the only angular point can exists (no angular points is also possible,
if no branches attach to the trunk from the left side). If there is an angular point, we cut
the strip $\Pi_1$ into two semi-strips by a segment $t=t_1$ through the
angular point. For the second ODE there is the similar strip $\Pi'_1$
which is divided by the angular point on the section $t=t'_1.$ For each
semi-strip $t\ge t_1$ and $t\le t_1$ we construct the Lyapunov
function with its levels. This is possible if the neighborhood is chosen thin
enough. After that we define some equimorphism on the boundary curve $L \to
L'$ and extend it to the strip as equimorphism by the same construction as
in the previous section. This equimorphism is extended till the closure of
the strip, hence it is defined on $\gamma_1$.

To the left from $\mathcal U$ there is the only $\mathcal V$, its
rightmost boundary curve $R$ has a nearest IC $\beta$ being globally stable on $\R.$
Any IC passing through $L$ at some $t_-$ intersects $R$ at an only point at some $t=t_+ >
t_-$ and the difference $t_+ - t_-$ is uniformly bounded for all such ICs
due to Assumption 4. This defines the correspondence between $L$ and $R$
and hence defines an uniformly continuous map $R\to R'$ to the related
curve of the second ODE. After that we construct the extension of the
equimorphism into the strip between $R$ and $\beta$ in the same way as
above.

For the upper semi-strip we do as in the preceding section:
we choose some uniformly continuous map $g$ on $L$ for $t\ge t_1$ defining the
map to the analogous left boundary curve $L'$ for the second ODE for $t\ge t_1'$. Then for
any point $m$ inside $\Pi_1$ we take the IC through $m$. This IC intersect
$L$ at some point $m_b$. Let $t(m)\ge t_1$ be the related time coordinate.
The map $g$ transform $m_b$ to the point $m'_b$
on $L'$ through which an IC of the second ODE passes. We choose on this IC
a point with the time coordinate $t=t' + t(m)-t_1.$

The second point is the extension of the equimorphism from two neighboring
semi-strips corresponding to two neighboring $u$-points of the equipped set into
a region between them being either of $V$-shape or of a semi-strip.
Let $\gamma_1$ and $\gamma_2$ be the related basic ICs for these
$u$-points. All IC between them belong to the same $s$-interval, that is
ICs leaving left semi-strip through its right boundary curve $r$ and ICs
leaving right semi-strip through its left boundary curve $l$ are asymptotic to
each other as $t\to \infty.$ Fix some IC $\gamma_*$ from this stable bunch (usually we
take as such IC that which go through an angular point where two boundary
curves of neighboring semi-strips intersect). For $t\ge T$ we choose some
thin uniform neighborhood of IC $\gamma_*$ such that this neighborhood
does not intersect both neighborhoods of $\gamma_1$ and $\gamma_2$. The
equimorphism constructed defines the map on the boundary curves and by ICs
it is extended to the boundary curves of the stable IC $\gamma_*$ and
after that does to its inner part.

\section{Autonomization}

In this section we introduce another invariant of uniform equivalency
which shows that a gradient-like ODE on $S^1$ is uniformly equivalent to
an asymptotically autonomous ODE. Let a gradient-like ODE be given.
Then a definite u-invariant corresponds to this ODE. Consider the circle
$S^1_0$ and related points $(u_1,\ldots,u_n,s_1,\ldots s_m)$. The numeration
in each group (``u'' and ``s'') is done in accordance with the orientation of $S^1$.
The section $S^1_0$ is
divided into $n$ intervals with points $u_i:$ $(u_1,u_2),$ $(u_2,u_3),\ldots,$
$(u_n,u_1)$. Each interval $(u_i,u_{i+1})$ contains some number $k_i$ $s$-points
arranged in the orientation order, $k_1+k_2+\cdots + k_n = m.$ Observe
that some $k_i$ may be zeroes.

Next we choose a smooth function $f_+$ on $S^1$ whose all zeroes are simple. In particular, zeroes at the points
$(u_1,\ldots,u_n)$ are all with positive derivatives and no other zeroes of this type exist.
Hence, more $n$ simple zeroes $s'_i$, $i= 1,\ldots,n,$ with negative derivatives should exist
at some points of each interval $(u_i,u_{i+1})$. The collection of points $(u_1,s'_1,\ldots,u_n,s'_n)$
makes up a complete set of zeroes for $f_+.$

Similarly, there is a smooth function $f_-$ with only simple zeroes,
and its zeroes with negative derivatives are just points $(s_1,\ldots,s_m).$
Then there are exactly $m$ simple zeroes $(u'_1,u'_{m})$ with positive derivatives at some
points on every interval $(s_j,s_{j+1}).$ Again, we assume these be all zeroes of $f_-.$
Functions $f_+, f_-$ define autonomous ODEs
at $\pm\infty.$ The scheme of the asymptotically autonomous ODE under construction is the following. Take an
annulus (cylinder) $K$ on the plane $\R^2$ with polar coordinates $(r,\varphi):$ $1\le r \le
2.$ We mark at the outer boundary $r=2$ zeroes of the function $f_+$
in accordance with the agreement that $\varphi/2\pi$ corresponds to the
coordinate on $S^1$. Similarly, we mark zeroes of $f_-$ at the inner
boundary $r=1$ with the same agreement. Now we join points
$(u_1,\ldots,u_n)$ on $r=2$ by the disjoint paths with points $(s_1,\ldots
s_m)$ on the boundary $r=1$ in accordance with the equipped set. Namely,
take the point $u_1$ of the equipped set. There is the unique IC with
unstable type exponential dichotomy on $\R_+.$ This IC belongs to a unique
bunch of ICs which possess the unstable type of dichotomy on $\R_-$. This
bunch defines two boundary ICs of ODE with the stable type of dichotomy on
$\R_-$, that is, two neighboring points $s_j,s_{j+1}$ of the equipped set.
On the interval in the circle $r=1$ between its related zeroes
$s_j,s_{j+1}$ with negative derivatives there is only one zero $u'_k$ with
the positive derivative. We joint point $u_1$ on $r=2$ with $u'_k$ on $r=1$ by a
path. In order to avoid intersections of the paths let us cut the annulus
along the path constructed. We get a curvilinear rectangle with two
boundary curves, $r=2$ (top) and $r=1$ (bottom), and two side boundaries
corresponding to the path. We assume that the orientation of top and
bottom boundaries corresponds to increasing $r$. Then $u_1$ is at the left
upper vertex of the rectangle and $u'_k$ is at the left lower vertex.

Take now $u_2$. If in the equipped set for ODE no $s$-points between $u_1$
and $u_2$ then both $u_1, u_2$ belong to the same
interval of $S^1_0$ defined by the $s$-points. This means, due to the
construction, that both ICs corresponding to $u_1,u_2$ belong to the same
unstable bunch. This means that as $t\to -\infty$ both these ICs are
asymptotically approach to each other. Then we need to join with a path
point $u_2$ on $r=2$ with the same point $u'_k$, this path should not
intersect the side boundaries, of course. This second path distinguishes a
curvilinear triangle in the rectangle and if we cut the rectangle along
this second path then we again obtain a rectangle but with the left upper
vertex $u_2.$ Suppose now that between $u_1$ and $u_2$ in the equipped set
there is at least one $s$-point, say $s_1$ (re-enumerating the set of $s$-points
one can always regard this). On the interval in $S^1_0$ defined by two
neighboring points $(u_1,u_2)$ we get the ordered set of $s$-points
$s_1,s_2,\ldots,s_{k_1}.$ The next $s$-point $s_{k_1+1}$ lies to the right from $u_2.$
Now we choose the interval on the bottom boundary $r=1$ between two zeroes
with negative derivatives of the function $f_-$ corresponding
$s_{k_1},s_{k_1+1}$. Due to the choice of $f_-$, between these two zeroes
there is a unique zero with the positive derivative $u'_2$. Joining by a path points
$u_2$ on $r=2$ with the point $u'_2$ on $r=1$ gives a sub-rectangle in the
rectangle with the side boundaries: the first path (left) and new path
(right). Cutting off the sub-rectangle we come to the new rectangle with
the left upper vertex $u_2$ instead of $u_1$ and $s_{k_1+1}$ being the
left lower vertex. To finish the first step of the induction we connect
within the sub-rectangle obtained points $s_1,s_2,\ldots,s_{k_1}$ by paths
with the unique point $s'_1$ on $r=2$ corresponding to a zero of $f_+$ with the negative
derivative on the interval $(u_1,u_2)$.

\begin{figure}[h]
\begin{center}
\includegraphics[width=15cm, keepaspectratio,
scale=4.0]{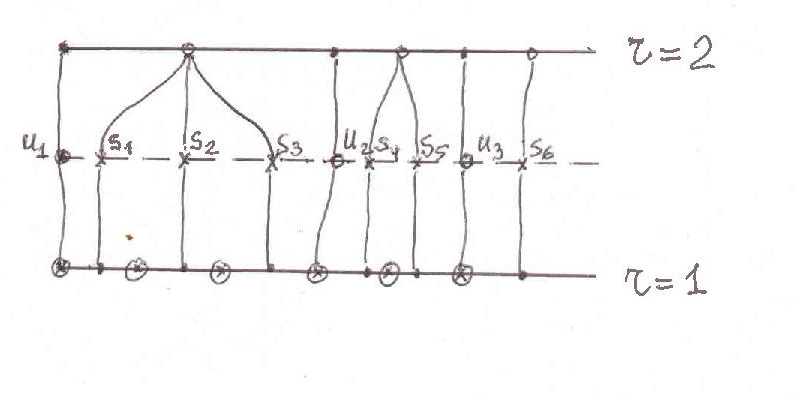}
\caption{Construction of asymptotically autonomous ODE}
\label{6}
\end{center}
\end{figure}

Thus we have made one step of the induction. Repeating this
procedure gives the annulus with all paths (see Fig.\ref{6}). This construction gives the
asymptotically autonomous gradient-like ODE with its
limiting ODEs defined by functions $f_+, f_-.$
The following theorem states this.
\begin{theorem}
For any gradient-like nonautonomous ODE $v$ on $S^1$ there is an asymptotically autonomous gradient-like
ODE being uniformly equivalent to $v.$
\end{theorem}
\proof To construct such an ODE we take on the cylinder $\R\times S^1$ the
foliations generated by $\dot x = f_-$ for $t< -T$ and the foliation for $\dot x =
f_+$ when $t>T.$ In the annulus $t= -T$ and $t=T$ we construct smooth
foliation as in the example above in accordance with behavior of paths
constructed. This gives the needed asymptotically autonomous ODE.
$\blacksquare$

\section{Almost periodic gradient-like equations on $S^1$}

In this section we apply the theory developed in previous sections to the classical
nonautonomous case -- almost periodic in time scalar
differential equations on the circle $S^1=\R/\Z$
 \begin{equation}\label{ap}
 \dot{x}=f(t,x),\; f(t,x+ 1)\equiv f(t,x).
 \end{equation}
We assume, as above, the equation be gradient-like one. The equation (\ref{ap}) is generated by the map $\hat{f}: \R \to
V^r(S^1)$, $r\ge 1,$ being almost periodic and therefore uniformly continuous.

Let us recall some notions from the theory of almost periodic functions $\cite{Levitan,Corduneanu,Zhikov}$.
A continuous map $h:\R \to B$ into Banach space $B$ is called {\em almost periodic}, if
for any $\varepsilon > 0$ there exists a relatively dense set $L(\eps)\subset
\R$ of $\eps$-almost periods such that for any $l\in L(\eps)$ the following
inequality holds
$$
\mathop{{sup}}\limits_{t\in \R}||h(t+l)-h(t)||< \eps.
$$
Recall that a subset $L\subset \R$ is relatively dense, if there is a positive
number $T \in \R$ such that any interval $(a,a+T),$ $a\in \R,$ contains at least one
number $l\in L.$

There is another definition of an almost periodic function introduced by Bochner.
It relies on the following theorem by Bochner \cite{Bochner}. For a
bounded continuous function $h: \R \to B$ consider the sequence of shifted functions $\{h(t+t_n)\}$.
\begin{definition}
Function $h$ is called normal, if for any sequence of its
shifts $\{h(t+t_n)\}$ there is a subsequence $\{h(t+t_{n_k})\}$ such that this
subsequence converges in the topology of uniform convergence on
$\R.$
\end{definition}

\begin{theorem}[Bochner]
A continuous function $h:\R \to B$ is almost periodic iff it is normal.
\end{theorem}
Trivial examples of gradient-like almost periodic equation on $S^1$ are
the following ones. The first is autonomous differential equation on $S^1$, $\dot x = f(x)$, when a $C^1$-smooth
1-periodic function has only simple zeroes (where $f' \ne 0$). The second is
when $f(t,x),$ $f(t,x+1)\equiv f(t,x),$ is $T$-periodic in $t$ and the related Poincar\'e map on $S^1$ in the period
$T$ is a diffeomorphism of $S^1$ with a rational rotation number. Then, as is known, all solutions are either
$nT$-periodic
or tend to $nT$-periodic solutions as $t\to \pm\infty.$ The integer $n\in\N$ is the same for all periodic solutions.
If all periodic
solutions have multipliers distinct of unity (i.e. the related periodic points of the Poincar\'e map are hyperbolic),
then this
periodic ODE is gradient-like one. A somewhat less trivial example is a small smooth almost periodic
perturbation of any of these two examples.

The second example shows that an almost periodic solution can go around $S^1$
infinitely many times when $t$ increases (decreases). For instance, any function of the form $x(t)=\alpha t +
u(t)\;(\mbox{\rm\;mod}\;1)$ with any $\alpha \in \R$ and arbitrary almost periodic
function $u$ is almost periodic on $S^1$. In fact, in this case we have a
particular example of an almost periodic mapping which takes its values in
some metric space $(X,d)$ with the metrics $d$. Recall that a continuous mapping $g: \R \to
X$ is almost periodic, if for any $\eps > 0$ there is a relatively dense
set $L(\eps)\subset \R$ such that for any $l\in L(\eps)$ the following
inequality holds
$$
d(g(t+l),g(t))\le \eps \;\mbox{for any\;} t\in \R.
$$
For our case $X=S^1= \R/\Z$ with the standard metrics $d(x,y)= \min\{|x-y|,
1-|x-y|\}$, $0\le x,y \le 1.$

Another equivalent definition of a.p. mapping $g: \R \to X$ is due to Bochner:
any sequence of shifts $g(t+\tau_n)$ is precompact in the topology of the
uniform convergence of mappings $\R \to X.$

For the case of $S^1$ a function $x(t)$ be an almost periodic function with values on $S^1$
if it is almost periodic as a function $\R\to \C$ understanding $S^1$
as being embedded into $\C$ as $|z|=1.$

In this section we prove the following theorem
\begin{theorem}
Suppose a differential equation $\dot{x}=f(t,x)$ on $S^1$ is gradient-like and almost periodic in $t$.
Then each its solution being unstable on $\R_+$ is almost periodic and possesses
exponential dichotomy of the unstable type on $\R$,
and each its solution being stable on $\R_-$ is almost periodic one and possesses
exponential dichotomy of the stable type on
$\R$. These a.p. solutions can be united into finitely many $n\in \N$ pairs -- one
stable and one unstable -- serving the boundary solutions of a closed strip.
The union of these strips is $\R\times S^1$.
Within such a strip solutions distinct of boundary ones tend to the boundary
stable almost periodic solution as $t\to \infty$
and to the boundary unstable a.p. solution as $t\to
-\infty$. This ODE is uniformly equivalent to the
autonomous scalar differential equation $\dot x = f(x)$ with a smooth function $f$ having
$2n$ simple zeroes on $S^1$.
\end{theorem}

To begin with we formulate several auxiliary statements.
\begin{lemma}
Suppose $\gamma$ be an integral curve possessing exponential dichotomy of the unstable type on $\R_+$.
Then there is an uniform neighborhood of $\gamma$
in $\R_+ \times S^1$ such that any IC intersecting this neighborhood belongs to the neighborhood only during a finite
time on $\R_+$. The boundary curves of this neighborhood are uniformly transversal to ICs intersecting them.
\end{lemma}
Such a neighborhood will be called an {\it isolating} neighborhood of such integral
curve. As such a neighborhood of the IC related to $\gamma$ the
neighborhood defined by a level of the Lyapunov function corresponding to
the solution $x_0(t)$, corresponding to $\gamma$, can be chosen (see above).
\begin{lemma}
Suppose some IC of the unstable type on $\R_+$ for ODE (\ref{ap}) is given.
There is $\delta>0$ such that if $||f(t,x)-g(t,x)||_{V^1(S^1)}<\delta,$ then the boundary
curves of an isolating neighborhood $U$ for this IC remain
uniformly transversal for IC of the equation $\dot x = g(t,x)$ and this
neighborhood contains the only IC $x_g(t)$ which stay wholly in $U$ for all $t\in \R_+.$
\end{lemma}
\begin{lemma}
IC $x_g(t)$ depends continuously on $||f(t,x)-g(t,x)||_{V^1(S^1)}$ at the ``point''
$f$.
\end{lemma}
Proofs of these lemmata can be obtained from standard results of the
theory of exponential dichotomy, see \cite{Demid}.

Consider some solution $x_0(t)$ of the differential equation
(\ref{ap}) which possesses exponential dichotomy of the unstable type on $\R_+$.
Since this equation is gradient-like one, such solution exists and, in virtue of Assumption 3, it also possesses
exponential dichotomy of the unstable type on $\R_-.$
This implies this solution to possess the exponential dichotomy of unstable type on the whole $\R.$ Choose some its isolating
neighborhood $U$ on $\R$ and consider any sequence $\tau_n \to \infty$ as $n\to \infty.$ We need to prove
that the sequence of shifted functions $x_0(t+\tau_{n})$ is precompact, i.e. a subsequence
$n_k \to \infty$ exists such that the subsequence $x_0(t+\tau_{n_k})$
converges in the topology of uniform convergence on $\R.$

Since ODE (\ref{ap}) is almost periodic, there exists some subsequence
$n_k \to \infty$ as $k\to \infty$ such that
the sequence of shifts $f(t+\tau_{n_k},x)$ converges in the topology of
uniform convergence on $\R$ in the space of continuous maps $\R \to V^r(S^1)$.
Thus, for any $\eps>0$ there is a $K(\eps) \in \N$ such that for any $k,k'> K(\eps)$
the inequality holds: $||f(t+\tau_{n_k},x)- f(t+\tau_{n_{k'}},x)||< \eps,$ where
the norm is taken in the space $V^r(S^1).$
The differential equation $\dot x = f(t+\tau_{n_k},x)$ evidently has the solution $x_0(t+\tau_{n_k}).$
If $U$ is an isolating neighborhood for IC corresponding to
$x_0(t)$, then its shift on $t_{n_k}$ gives the isolating neighborhood of shifted
solution. For $\eps$ small enough ODEs with r.h.s.
$f(t+\tau_{n_k},x)$ and $f(t+\tau_{n_k'},x)$ are close enough and since
$f$ is rough, then all shifted equations are also rough with the same radius of roughness.
In particular, for $\eps$ small enough an equimorphism realizing the uniform equivalence
of the foliation onto ICs for ODE $\dot x = f(t+\tau_{n_k},x)$ with the foliation for ODE
$\dot x = f(t+\tau_{n_{k'}},x)$ is close to the identity map $id_{\R\times S^1}$. This implies
that in the isolating neighborhood of the solution
$x_0(t+\tau_{n_k})$ can lie wholly only one solution of ODE with r.h.s. $f(t+\tau_{n_{k'}},x)$
and this solution has the exponential dichotomy of unstable type on $\R$.
Thus, we get a unique solution of the ODE $\dot x = f(t+\tau_{n_{k'}},x)$
with exponential dichotomy of unstable type on $\R.$ Shifting it back in time allows us to get
a solution $x_1(t)$ of the same type for the initial ODE (\ref{ap}). We can do the
same procedure for $x_0(t)$ changing the sequence $\tau_n$ and sub-sequences
$n_k$. In such a way, we shall find several solutions with the same type of
dichotomy of the initial ODE which are connected by the limiting procedure
with $x_0(t)$. Since there are finitely many solutions with
the dichotomy of unstable type on $\R$, we break all such solutions into
several groups. It is clear that if, using this procedure, we get a solution
$x_1(t)$ from the solution $x_0(t)$ and next do it from $x_1(t)$ to $x_2(t)$,
then combining a proper sequence of shifts we get $x_2(t)$ from $x_0(t),$ as well.
Thus, solutions from different groups cannot be obtained from each other by the limit.
But all solutions from the same group can be obtained from each other by the given procedure. In
particular, this means that starting from solution $x_0(t)$ we can get
$x_0(t)$ itself using some sequences of shifts. In a sense, one group of
solutions corresponds to one periodic orbit and its shifts on the periods
for the periodic ODE. Observe, that solutions from one group are separated from each other:
these solutions are separated by their uniform isolating neighborhoods.

Thus we proved that shifts of a given solutions
$x_0(t)$ with the exponential dichotomy of unstable type on $\R$ is
precompact: its limit sets consists of solutions of one group. Therefore,
this solution is almost periodic.

The same considerations work if we start from a solution of ODE (\ref{ap})
with the exponential dichotomy of stable type on $\R_-$. Then we conclude,
due to gradient-likeness of (\ref{ap}), that this solution possesses by
the exponential dichotomy of stable type on $\R$. After that we again
break solutions of such type into groups, and so forth. Solutions from
different groups of stability cannot coinside and therefore all extended phase
manifold $\R\times S^1$ is divided into several strips whose boundary ICs
are one stable and one unstable ones.
$\blacksquare$

As an immediate application of this theorem, let us consider an almost periodic perturbation of a periodic
ODE on $S^1$. Concerning the unperturbed ODE we assume that its Poincar\'e map in the period is rough, that is it has
a rational rotation number and all its periodic points are hyperbolic.
Then this periodic equation is gradient-like and rough. Under a sufficiently small a.p.
perturbation near each hyperbolic periodic IC in $\R\times S^1$ a unique
a.p. IC arises. The type of exponential dichotomy for the a.p. solutions on $\mathbb R$
will be the same as was for the related periodic solution of unperturbed
equation. All other solutions will tend to a.p. solutions as $t\to \pm \infty.$
Thus we get
\begin{corollary}
A sufficiently small a.p. perturbation of a periodic ODE with a rough
Poincar\'e map in the period is uniformly equivalent to the autonomous ODE
with the same number of its simple zeroes as for the number of periodic orbits
for the Poincar\'e map.
\end{corollary}

It is worth remarking that if an almost periodic ODE fails to be rough
then it may not have any almost periodic solutions at all. For instance,
a periodic ODE on $S^1$ having an irrational rotation number and realizing
the Denjoy case (for $C^1$-smooth $f$) has not any almost periodic solutions.
The first such example for an almost periodic case was constructed in \cite{Opial}
(see its more detailed consideration in \cite{Ler_rem}).

\section{Acknowledgement}

L.M.L. thanks for a support the Russian Science Foundation (grant 14-41-00044), as well as
the Russian Foundation of Basic Research (grants 16-51-10005-КО and 16-01-00324). A part of this
research was supported by the Ministry of Education and Science of Russian Federation (project
1.3287.2017, target part).

\end{document}